\documentclass[a4paper,11pt,twoside,reqno]{amsart}

\usepackage[utf8]{inputenc}
\usepackage[plainpages=false,pdfpagelabels=true]{hyperref}
\usepackage{amssymb,amsthm,slashed}
\usepackage[margin=1in]{geometry}

\newtheorem{Satz}{Theorem}[section]
\newtheorem{Prop}[Satz]{Proposition}
\newtheorem{Lem}[Satz]{Lemma}

\newtheorem{Cor}[Satz]{Corollary}
\theoremstyle{definition}
\newtheorem{Dfn}[Satz]{Definition}
\newtheorem{Bem}[Satz]{Remark}
\newcommand{\sff}{\mathrm{I\!I}}

\newcommand{\p}{\slashed{\partial}}
\newcommand{\D}{\slashed{D}}

\newcommand{\cL}{\mathcal{L}}
\allowdisplaybreaks[1]

\renewcommand{\epsilon}{\varepsilon}

\newcommand{\R}{\ensuremath{\mathbb{R}}}

\newcommand{\C}{\ensuremath{\mathbb{C}}}

\numberwithin{equation}{section}

\title{On conservation laws for the supersymmetric sigma model}
\author{Volker Branding}
\date{\today}
\address{University of Vienna, Faculty of Mathematics\\
Oskar-Morgenstern-Platz 1, 1090 Vienna, Austria\\}
\email{volker.branding@univie.ac.at}
\subjclass[2010]{58E20, 53C27, 70S10}
\keywords{supersymmetric nonlinear sigma model; Dirac-harmonic maps; conservation laws}
\begin{document}

\begin{abstract}
We derive conservation laws for Dirac-harmonic maps and their extensions
to manifolds that have isometries, where we mostly focus on the spherical case.
In addition, we discuss several geometric and analytic applications of the 
latter.
\end{abstract} 

\maketitle

\section{Introduction and Results}
Symmetries have always been a driving principle in both mathematics and physics.
This statement manifests itself as \emph{Noether's Theorem}, namely that every
continuous symmetry of a system leads to a conservation law.

When constructing a physical model to describe elementary particles
one considers an energy functional together with a certain amount of symmetries 
that leave it invariant. These symmetries can both be discrete and continuous.
The energy functionals considered in quantum field theory are formulated in terms
of objects from differential geometry. Consequently, their invariances and the 
corresponding conversation laws also allow for a geometric interpretation.

When considering an energy functional that involves a map between two manifolds,
the invariance under diffeomorphisms on the domain gives rise to the energy-momentum tensor,
which is conserved for a critical point.
Moreover, symmetries on the target lead to a different conserved quantity,
which is called Noether current in the physics literature.

Throughout this article we will study an action functional that is motivated
from the supersymmetric nonlinear sigma model from quantum field theory \cite{MR626710},
see also \cite{MR1858211}.
From a geometric point of view this energy functional consists of the energy for harmonic maps
coupled to spinor fields in a nonlinear fashion. 
For a recent survey on harmonic maps we refer to \cite{MR2389639}, for an introduction
to spin geometry see \cite{MR1031992}.
The geometric study of the supersymmetric sigma model was initiated in \cite{MR2262709,MR2176464},
where the notion of \emph{Dirac-harmonic maps} was introduced. This notion was extended later on to
include an additional curvature term \cite{MR2370260, MR3333092}, a two-form potential \cite{MR3305429}
and a connection with metric torsion on the target \cite{MR3493217}.
Currently, many analytic and geometric aspects of Dirac-harmonic maps and their extensions
are well-understood, like the regularity of weak solutions \cite{MR2544729,MR3333092}.
However, apart from an existence result for uncoupled Dirac-harmonic maps 
\cite{MR3070562}, a general existence result is still not available. 

For a general treatment of harmonic maps and conservation laws we refer to the book \cite{MR1913803}.
For a supergeometric study of harmonic maps, see \cite{MR2316716}.
In this article we focus on the derivation of conservation laws
for critical points of the supersymmetric nonlinear sigma model
to targets with symmetries.

It is well known that both nonlinear Dirac equations on surfaces 
and harmonic maps from surfaces to spheres have a natural connection to CMC surfaces.
The critical points of the supersymmetric nonlinear sigma model interpolate
between these equations and we discuss a geometric interpretation of the combined system.

This article is organized as follows: 
In Section 2 we recall the mathematical setup that we use to study Dirac-harmonic maps and Dirac-harmonic
maps with curvature term. In the third section we consider the case of a spherical target, derive
a conservation law and give several geometric and analytic applications.
Section 4 is then devoted to a target with isometries, where we again derive a conservation law
for Dirac-harmonic maps and Dirac-harmonic maps with curvature term.

\section{The nonlinear supersymmetric sigma model as a geometric variational problem}
Let us describe the mathematical setup used in this article.
Let \((M,h)\) be a closed Riemannian spin manifold with spinor bundle \(\Sigma M\).
We fix both a spin structure and a metric connection \(\nabla^{\Sigma M}\) on \(\Sigma M\). 
Moreover, we fix a hermitian scalar product on \(\Sigma M\) which we denote by \(\langle\cdot,\cdot\rangle_{\Sigma M}\).
On the spinor bundle there is the Clifford multiplication of spinors with tangent vectors, which is skew-symmetric
\begin{align*}
\langle\psi,X\cdot\eta\rangle_{\Sigma M}=-\langle X\cdot\psi,\eta\rangle_{\Sigma M}
\end{align*}
for all \(X\in TM,\eta,\psi\in\Gamma(\Sigma M)\).
Moreover, we have the Clifford relations, that is
\begin{align*}
X\cdot Y\cdot\psi+Y\cdot X\cdot\psi=-2h(X,Y)\psi
\end{align*}
for all \(X,Y\in TM\).
The Dirac operator \(\p\) acting on sections of \(\Sigma M\) is defined as
\begin{align*}
\p:=e_\alpha\cdot\nabla^{\Sigma M}_{e_\alpha}, 
\end{align*}
where \(e_\alpha,\alpha=1,\ldots,\dim M\) is a basis of \(TM\).
Throughout this article we will make use of the summation convention, that is we sum over repeated indices.
The Dirac operator \(\p\) is elliptic and self-adjoint with respect to the \(L^2\)-norm.

We will mostly consider a two-dimensional domain \(M\), in this case the spinor bundle
splits as \(\Sigma M=\Sigma^+M\oplus\Sigma^-M\), where we call \(\Sigma^+M\) the bundle
of \emph{positive spinors} and \(\Sigma^-M\) the bundle of \emph{negative spinors}.
We will make use of the complex volume element \(\omega_\C\), which is defined by
\begin{align*}
\omega_\C:=ie_1\cdot e_2.
\end{align*}
In order to project to the subbundles \(\Sigma^\pm M\) we make use of the projector
\begin{align}
\label{projector-chiralities}
P_\pm:=\frac{1}{2}(1\pm\omega_\C).
\end{align}

In addition, let \((N,g)\) be another Riemannian manifold.
Consider a map \(\phi\colon M\to N\), which we use to pullback the tangent bundle \(TN\) to \(M\).
We form the twisted bundle \(\Sigma M\otimes\phi^\ast TN\), 
sections in this bundle will be called \emph{vector spinors}. We will denote the connection on \(\Sigma M\otimes\phi^\ast TN\)
by \(\tilde{\nabla}\). This leads to the twisted Dirac operator acting on vector spinors, which is given by
\begin{align*}
\D:=e_\alpha\cdot\tilde{\nabla}_{e_\alpha}.
\end{align*}
This twisted Dirac operator is also elliptic and self-adjoint with respect to the \(L^2\)-norm.

If we choose local coordinates we will use Greek indices for coordinates on the domain \(M\)
and Latin indices for coordinates on the target \(N\).
Whenever clear from the context we will use \(\langle\cdot,\cdot\rangle\) for a generic scalar product
without referring to the actual bundle. 

\subsection{Dirac-harmonic maps and extensions}
In this section we recall the action functional that we
will mostly investigate in this article
\begin{align}
\label{energy-kappa}
E_\kappa(\phi,\psi)=\int_M(|d\phi|^2+\langle\psi,\D\psi\rangle+\kappa\langle R^N(\psi,\psi)\psi,\psi\rangle)dM.
\end{align}
Here, \(\kappa\) is a real-valued parameter. The first term is the usual harmonic energy for a map
between two Riemannian manifolds, in the second term the scalar product is taken on the bundle \(\Sigma M\otimes\phi^\ast TN\).
In the third term the spinors are contracted as
\[
\langle R^N(\psi,\psi)\psi,\psi\rangle_{\Sigma M\otimes\phi^\ast TN}=R_{ijkl}\langle\psi^i,\psi^k\rangle_{\Sigma M}\langle\psi^j,\psi^l\rangle_{\Sigma M},
\]
which ensures that the action is real-valued. Here \(R_{ijkl}\) are the components of the Riemann curvature tensor on \(N\).
The critical points of \eqref{energy-kappa} are given by (see \cite[Proposition 2.1]{MR3333092})
\begin{align}
\label{kappa-critical-phi}\tau(\phi)=&\frac{1}{2}R^N(e_\alpha\cdot\psi,\psi)d\phi(e_\alpha)+\frac{\kappa}{2}\langle(\nabla R^N)^\sharp(\psi,\psi)\psi,\psi\rangle, \\
\label{kappa-critical-psi}\D\psi=&-2\kappa R^N(\psi,\psi)\psi.
\end{align}
Note that the energy functional \eqref{energy-kappa} and its critical points \eqref{kappa-critical-phi}, \eqref{kappa-critical-psi} 
interpolate between the energy functionals for Dirac-harmonic maps (\(\kappa=0\)) and Dirac-harmonic maps with curvature term (\(\kappa=-\frac{1}{6}\)).

We call solutions of the system \eqref{kappa-critical-phi}, \eqref{kappa-critical-psi}
\(\kappa\)-\emph{Dirac-harmonic maps}.

We want to point out that the energy functional we are considering here 
is slightly different compared to the physics literature since we do not use Grassmann-valued spinors.
However, this leads to the advantage that we can employ well-established methods from geometric analysis to
study \eqref{energy-kappa} and its critical points.

In the following we will sometimes need the extrinsic version of \eqref{kappa-critical-phi}, \eqref{kappa-critical-psi}.
To this end we apply the Nash embedding theorem to isometrically embed \(N\) into some \(\R^q\).
We will denote the second fundamental form of the embedding by \(\sff\).

The extrinsic version of \eqref{kappa-critical-phi}, \eqref{kappa-critical-psi} is given by the system
\begin{align}
\label{phirq}-\Delta\phi=&\sff(d\phi,d\phi)+P(\sff(e_\alpha\cdot\psi,d\phi(e_\alpha)),\psi)-\kappa G(\psi), \\
\label{psirq}\p\psi=&\sff(d\phi(e_\alpha),e_\alpha\cdot\psi)+\kappa F(\psi,\psi)\psi
\end{align}
with the terms
\begin{align*}
G(\psi)&=-\frac{1}{2}(\langle(\nabla\sff)(\partial_{y^i},\partial_{y^k}),\sff(\partial_{y^j},\partial_{y^l})\rangle-
\langle(\nabla\sff)(\partial_{y^i},\partial_{y^l}),\sff(\partial_{y^j},\partial_{y^k})\rangle)\langle\psi^i,\psi^k\rangle\langle\psi^j,\psi^l\rangle, \\
F(\psi,\psi)\psi&=-2(P(\sff(\partial_{y^k},\partial_{y^j}),\partial_{y^l})-P(\sff(\partial_{y^l},\partial_{y^j}),\partial_{y^k}))\langle\psi^j,\psi^l\rangle\psi^k,
\end{align*}
where now \(\phi\colon M\to\R^q\) and \(\psi\in\Gamma(\Sigma M\otimes\R^q)\).
In addition, \(P\) denotes the shape operator that is defined via
\begin{align*}
\langle P(\xi,X),Y\rangle_{TN}=\langle\sff(X,Y),\xi\rangle_{\R^q}
\end{align*}
for all \(X,Y\in TN\) and \(\xi\in T^\perp N\).

The extrinsic version allows us to consider a weak formulation of \eqref{kappa-critical-phi}, \eqref{kappa-critical-psi}.
To this end, we define the following space 
\begin{align*}
\chi(M,N):=\{(\phi,\psi)\in W^{1,2}(M,N)\times W^{1,\frac{4}{3}}(M,\Sigma M\otimes\phi^\ast TN) 
\text{ with } \eqref{kappa-critical-phi} \text{ and } \eqref{kappa-critical-psi} \text{ a.e.}\}.
\end{align*}
\begin{Dfn}
A pair \((\phi,\psi)\in\chi(M,N)\) is called \emph{weak \(\kappa\)-Dirac-harmonic map} from \(M\) to \(N\) if and only if the pair \((\phi,\psi)\) solves
\eqref{kappa-critical-phi}, \eqref{kappa-critical-psi} in a distributional sense.
\end{Dfn}

\subsection{Spinorial symmetries}
Before we discuss how isometries on the target lead to conservation laws 
we briefly discuss some symmetries arising in the context of the spinors.
To this end we recall the following
\begin{Lem}
\label{lemma-complex-volume-element}
The complex volume element \(\omega_\C\) satisfies the following algebraic relations
\begin{enumerate}
 \item \(\omega_\C^2=1\).
 \item \(X\cdot\omega_\C=-\omega_\C\cdot X\) for a two-dimensional manifold and \(X\in TM\).
 \item \(\langle\omega_\C\cdot\psi,\eta\rangle_{\Sigma M}=-\langle\omega_\C\cdot\psi,\eta\rangle_{\Sigma M}\)
  for all \(\psi,\eta\in\Sigma M\).
\end{enumerate}
\end{Lem}

\begin{Prop}
The energy functional \eqref{energy-kappa} and its critical points \eqref{kappa-critical-phi}, \eqref{kappa-critical-psi},
are invariant under the symmetries
\begin{enumerate}
 \item \(\psi\to e^{i\alpha}\psi\) with \(\alpha\in\R\)
 \item \(\psi\to \omega_\C\cdot\psi\)
\end{enumerate}
\end{Prop}

\begin{proof}
Note that the complex volume element \(\omega_\C\) is parallel 
and thus by Lemma \ref{lemma-complex-volume-element} we find
\[
\D(\omega_\C\cdot\psi)=-\omega_\C\cdot\D\psi.
\]
Consequently, we obtain
\[
\langle\omega_\C\cdot\psi,\D(\omega_\C\cdot\psi)\rangle_{\Sigma M\otimes\phi^\ast TN}=-\langle\omega_\C\cdot\psi,\omega_\C\cdot\D\psi\rangle_{\Sigma M\otimes\phi^\ast TN}
=\langle\psi,\D\psi\rangle_{\Sigma M\otimes\phi^\ast TN}
\]
and
\begin{align*}
\langle\omega_\C\cdot\psi^k,e_\alpha\cdot\omega_\C\cdot\psi^k\rangle_{\Sigma M}
=-\langle\omega_\C\cdot\psi^k,\omega_\C\cdot e_\alpha\cdot\psi^k\rangle_{\Sigma M}=\langle\psi^k,e_\alpha\cdot\psi^k\rangle_{\Sigma M}, \\
\langle\omega_\C\cdot\psi^i,\omega_\C\cdot\psi^k\rangle_{\Sigma M}\langle\omega_\C\cdot\psi^j,\omega_\C\cdot\psi^l\rangle_{\Sigma M}
=\langle\psi^i,\psi^k\rangle_{\Sigma M}\langle\psi^j,\psi^l\rangle_{\Sigma M}
\end{align*}
proving the claim.
\end{proof}

Note that we also have discrete symmetries in the term
\[
R_{ijkl}\langle\psi^i,\psi^k\rangle_{\Sigma M}\langle\psi^j,\psi^l\rangle_{\Sigma M}.
\]

\begin{Bem}
In the physics literature the complex volume element \(\omega_\C\) is usually denoted by \(\gamma^5\).
\end{Bem}

\section{The case of a spherical target}
In this section we consider the system \eqref{kappa-critical-phi}, \eqref{kappa-critical-psi}
in the case of a spherical target.
On the one hand this particular case is attractive since the huge symmetry of the sphere easily leads to 
a conservation law. Moreover, in the physics literature nonlinear sigma-models are mostly considered 
having a spherical target.

For \(N=S^n\subset\R^{n+1}\) with the round metric the Euler-Lagrange equations read
\begin{align}
\label{kappa-sphere-phi}-\Delta\phi^i&=|d\phi|^2\phi^i+\operatorname{Re}\langle\psi^i,e_\alpha\cdot\psi^j\rangle_{\Sigma M}\phi^j_\alpha,\qquad i=1,\ldots,n+1, \\
\label{kappa-sphere-psi}\p\psi^i&=-\phi^j_\alpha e_\alpha\cdot\psi^j\otimes\phi^i-2\kappa(|\psi|_{\Sigma M}^2\psi^i-\langle\psi^j,\psi^i\rangle_{\Sigma M}\psi^j),\qquad i=1,\ldots,n+1, 
\end{align}
where we use the notation \(\phi^j_\alpha:=\frac{\partial\phi^j}{\partial x^\alpha}\).

\begin{Bem}
Suppose we have a smooth solution of the system \eqref{kappa-sphere-phi}, \eqref{kappa-sphere-psi}.
Using that \(|d\phi|^2=-\langle\Delta\phi,\phi\rangle\) for maps taking values in \(S^n\subset\R^{n+1}\) we find
\begin{align*}
\operatorname{Re}\langle\psi^i,e_\alpha\cdot\psi^j\rangle_{\Sigma M}\phi^j_\alpha\phi^i=0.
\end{align*}
If we think of the summation over \(i,j\) as taking the trace of
an endomorphism, then we may expect that the endomorphism itself contains
some interesting information.
\end{Bem}

In the following we show how the existence of isometries on the sphere 
leads to a conserved quantity. Becoming more technical, let us recall the following facts:
\begin{Dfn}
A vector field \(X\) is called \emph{Killing vector field} on \((N,g)\) if
\[
\mathcal{L}_Xg=0,
\]
where \(\cL\) represents the Lie-derivative of the metric. In terms of local coordinates we have
\[
0=(\mathcal{L}_Xg)_{ij}=g_{jk}\nabla_{\partial_{y^i}}X^k+g_{ik}\nabla_{\partial_{y^j}}X^k.
\]
\end{Dfn}

The group \(SO(n+1)\) acts isometrically on \(S^n\).
The set of Killing vector fields on \(S^n\) can be identified with the Lie algebra \(\mathfrak{so}(n+1)\) of \(SO(n+1)\).
In addition, \(\mathfrak{so}(n+1)\) can be represented as \((n+1)\times (n+1)\) skew-symmetric real-valued matrices.
For simplicity we will assume that these matrices have only entries of \(+1,0,-1\).

We will determine a conserved current in the case that we have a weak solution of \eqref{kappa-sphere-phi}, \eqref{kappa-sphere-psi},
where we follow the ideas from \cite{MR1085633} for harmonic maps. This method has the advantage of leading
to the result in a rather straightforward way.

\begin{Prop}
Let \((\phi,\psi)\colon M\to S^n\) be a weak \(\kappa\)-Dirac-harmonic map.
Then the following conservation law holds
\begin{equation}
\int_M\big(\operatorname{Re}\langle e_\alpha\cdot\psi^i,\psi^m\rangle_{\Sigma M}+(\phi^i_\alpha\phi^m-\phi^i\phi^m_\alpha)\big)(\nabla_{e_\alpha}\eta) dM=0
\end{equation}
for all \(\eta\in C^{\infty}(M)\).
\end{Prop}

\begin{proof}
Let \(X(\phi)\) be a Killing vector field on \(S^n\) and \(\eta\in C^\infty(M)\).
Testing \eqref{kappa-sphere-phi} with \(X^k(\phi)\eta\) we obtain
\begin{align*}
-\int_M\Delta\phi^kX^k(\phi)\eta dM=\int_M|d\phi|^2\underbrace{\phi^kX^k(\phi)}_{=0}\eta dM
+\int_M\operatorname{Re}\langle\psi^k,e_\alpha\cdot\psi^j\rangle_{\Sigma M}\phi^j_\alpha X^k(\phi)\eta dM.
\end{align*}
Note that the first term on the right hand side vanishes since \(X\perp\phi\).
In addition, we calculate
\begin{align*}
-\int_M\Delta\phi^kX^k(\phi)\eta dM
=\int_M\underbrace{\nabla\phi^k\nabla X^k(\phi)}_{=0} dM
+\int_M\nabla\phi^kX^k(\phi)\nabla\eta dM,
\end{align*}
where the second terms vanishes since \(X(\phi)\) is a Killing vector field.
Since Killing vector fields on the sphere can be identified with antisymmetric matrices, we 
find
\begin{align*}
\int_M\operatorname{Re}\langle\psi^k,e_\alpha\cdot\psi^j\rangle_{\Sigma M}\phi^j_\alpha X^k_{im}(\phi)\eta dM
=&\int_M(\operatorname{Re}\langle\psi^i,e_\alpha\cdot\psi^j\rangle_{\Sigma M}\phi^j_\alpha\phi^m \\
&-\operatorname{Re}\langle\psi^m,e_\alpha\cdot\psi^j\rangle_{\Sigma M}\phi^j_\alpha\phi^i)\eta dM \\
=&\int_M(\operatorname{Re}\langle\psi^i,\p\psi^m\rangle_{\Sigma M}
-\operatorname{Re}\langle\psi^m,\p\psi^i\rangle_{\Sigma M})\eta dM \\
=&\int_M(\operatorname{Re}(\underbrace{\langle\psi^i,\p\psi^m\rangle_{\Sigma M}-\langle\p\psi^m,\psi^i\rangle_{\Sigma M}}_{\operatorname{Im}\langle\psi^i,\p\psi^m\rangle_{\Sigma M}})\eta dM \\
&+\int_M\operatorname{Re}\langle\psi^i,e_\alpha\cdot\psi^m\rangle_{\Sigma M}(\nabla_{e_\alpha}\eta) dM \\
=&\int_M\operatorname{Re}\langle\psi^i,e_\alpha\cdot\psi^m\rangle_{\Sigma M}(\nabla_{e_\alpha}\eta) dM,
\end{align*}
which completes the proof.
\end{proof}

We can check by a direct calculation that given a smooth \(\kappa\)-Dirac-harmonic map
we obtain a vector field that is divergence free.

\begin{Lem}
Let \((\phi,\psi)\colon M\to S^n\) be a smooth \(\kappa\)-Dirac-harmonic map. Then the vector field
\begin{align}
\label{noether-current-kappa-sphere}
J_\alpha^{im}:=\operatorname{Re}\langle e_\alpha\cdot\psi^i,\psi^m\rangle_{\Sigma M}+(\phi^i_\alpha\phi^m-\phi^i\phi^m_\alpha)
\end{align}
is divergence free.
\end{Lem}
\begin{proof}
We calculate
\begin{align*}
\nabla_{e_\alpha}J_\alpha^{im}=&\operatorname{Re}\langle\p\psi^i,\psi^m\rangle_{\Sigma M}-\operatorname{Re}\langle\psi^i,\p\psi^m\rangle_{\Sigma M}
+\Delta\phi^i\phi^m-\Delta\phi^m\phi^i.
\end{align*}
Moreover, we find
\begin{align*}
\operatorname{Re}\langle\p\psi^i,\psi^m\rangle_{\Sigma M}
-\operatorname{Re}\langle\psi^i,\p\psi^m\rangle_{\Sigma M}
=&2\kappa\operatorname{Re}\big(-|\psi|_{\Sigma M}^2\langle\psi^i,\psi^m\rangle_{\Sigma M}
+\langle\psi^j,\psi^i\rangle_{\Sigma M}\langle\psi^j,\psi^m\rangle_{\Sigma M} \\
&+|\psi|_{\Sigma M}^2\langle\psi^i,\psi^m\rangle_{\Sigma M}
-\langle\psi^j,\psi^m\rangle_{\Sigma M}\langle\psi^j,\psi^i\rangle_{\Sigma M}\big) \\
&-\operatorname{Re}\langle e_\alpha\cdot\psi^j,\psi^m\rangle_{\Sigma M}\phi^j_\alpha\phi^i
+\operatorname{Re}\langle\psi^i,e_\alpha\cdot\psi^j\rangle_{\Sigma M}\phi^j_\alpha\phi^m \\
=&4\kappa\operatorname{Re}\big(\operatorname{Im}(\langle\psi^m,\psi^j\rangle_{\Sigma M}\langle\psi^i,\psi^j\rangle_{\Sigma M}\big) \\
&-\operatorname{Re}\langle e_\alpha\cdot\psi^j,\psi^m\rangle_{\Sigma M}\phi^j_\alpha\phi^i
+\operatorname{Re}\langle\psi^i,e_\alpha\cdot\psi^j\rangle_{\Sigma M}\phi^j_\alpha\phi^m \\
=&-\operatorname{Re}\langle e_\alpha\cdot\psi^j,\psi^m\rangle_{\Sigma M}\phi^j_\alpha\phi^i
+\operatorname{Re}\langle\psi^i,e_\alpha\cdot\psi^j\rangle_{\Sigma M}\phi^j_\alpha\phi^m.
\end{align*}
On the other hand we have
\[
\Delta\phi^i\phi^m-\Delta\phi^m\phi^i=-\operatorname{Re}\langle\psi^i,e_\alpha\cdot\psi^j\rangle_{\Sigma M}\phi^j_\alpha\phi^m
+\operatorname{Re}\langle\psi^m,e_\alpha\cdot\psi^j\rangle_{\Sigma M}\phi^j_\alpha\phi^i
\]
such that
\begin{align*}
\nabla_{e_\alpha}J_\alpha^{im}=&\phi_\alpha^j\phi^i\operatorname{Re}(\langle\psi^m,e_\alpha\cdot\psi^j\rangle_{\Sigma M}-\langle e_\alpha\cdot\psi^j\psi^m\rangle_{\Sigma M}) \\
=&2\phi_\alpha^j\phi^i\operatorname{Re}(\operatorname{Im}(\langle\psi^m,e_\alpha\cdot\psi^j\rangle_{\Sigma M}))\\
=&0,
\end{align*}
yielding the claim.
\end{proof}

Following the terminology used in the physics literature we call the vector field \(J\)
\emph{Noether current}. It is obvious that \(J\) is unique up to multiplication with an overall 
constant and the addition of a parallel vector field.

\begin{Bem}
The term \(\langle\psi^i,e_\alpha\cdot\psi^m\rangle_{\Sigma M}\) takes the form
of a Killing vector field associated to a Killing spinor.
More precisely, a Killing spinor \(\psi\) is a section of \(\Sigma M\) that satisfies
\[
\nabla_X\psi=\alpha X\cdot\psi,
\]
where \(\alpha\) is a non-vanishing complex number.
To a given Killing spinor we can associate a vector field \(V\) via the Riemannian metric
\[
h(V,X):=i\langle\psi,X\cdot\psi\rangle_{\Sigma M}.
\]
Not many Riemannian manifolds allow the existence of Killing spinors \cite{MR1224089}, however 
these always exist on spheres. Consequently, it is not surprising that a term having the form of
a Killing vector field appears in the Noether current for \(\kappa\)-Dirac-harmonic maps to spheres.
\end{Bem}

From now on we assume that \(M\) is two-dimensional and by \(\Omega\) we denote
a connected domain in \(M\). We denote the local coordinates on \(\Omega\) by \(x,y\)
and its tangent vectors by \(\partial_x,\partial_y\).

\begin{Bem}
Suppose that \(\phi\colon\Omega\to S^n\) is a smooth harmonic map.
In this case the Noether current reads
\[
J_\alpha^{im}:=\phi^i_\alpha\phi^m-\phi^i\phi^m_\alpha.
\]
By a direct calculation it follows that the Noether current satisfies the following algebraic relation
\[
\partial_x J_y-\partial_y J_x-2[J_x,J_y]=0,
\]
which can be thought of as a vanishing curvature condition if we think of \(J\) as 
the connection one-form on the bundle \(\phi^\ast TN\otimes\mathfrak{so}(n+1)\).
This fact relates the theory of harmonic maps to spheres to the world of integrable systems \cite{MR1264181}.
\end{Bem}

In the following we discuss if a similar structure also holds for \(\kappa\)-Dirac-harmonic maps
to spheres.

\begin{Lem}
The vector field \(J_\alpha^{im}\) defined in \eqref{noether-current-kappa-sphere} satisfies the following algebraic relation
\begin{align}
\label{noether-current-identity-general}
\partial_x J_y^{im}&-\partial_y J_x^{im}-2(J^{ij}_x J^{jm}_y-J^{ij}_y J^{jm}_x) \\
\nonumber=&-\operatorname{Re}\langle (\partial_x\cdot\nabla_{\partial_y}-\partial_y\cdot\nabla_{\partial_x})\psi^i,\psi^m\rangle_{\Sigma M}
+\operatorname{Re}\langle\psi^i,(\partial_x\cdot\nabla_{\partial_y}-\partial_y\cdot\nabla_{\partial_x})\psi^m\rangle_{\Sigma M} \\
\nonumber&-2\big(\operatorname{Re}\langle\psi^i,\partial_x\cdot\psi^j\rangle_{\Sigma M}\operatorname{Re}\langle\psi^j,\partial_y\cdot\psi^m\rangle_{\Sigma M}
-\operatorname{Re}\langle\psi^i,\partial_y\cdot\psi^j\rangle_{\Sigma M}\operatorname{Re}\langle\psi^j,\partial_x\cdot\psi^m\rangle_{\Sigma M} \\
\nonumber&-\operatorname{Re}\langle\psi^i,\partial_x\cdot\psi^j\rangle_{\Sigma M}\phi^j_y\phi^m
+\operatorname{Re}\langle\psi^j,\partial_y\cdot\psi^m\rangle_{\Sigma M}\phi^i\phi^j_x \\
\nonumber&+\operatorname{Re}\langle\psi^i,\partial_y\cdot\psi^j\rangle_{\Sigma M}\phi^j_x\phi^m
-\operatorname{Re}\langle\psi^j,\partial_x\cdot\psi^m\rangle_{\Sigma M}\phi^i\phi^j_y\big).
\end{align}
\end{Lem}

\begin{proof}
By a direct computation we find
\begin{align*}
\partial_x J_y^{im}-\partial_y J_x^{im}=&
-\operatorname{Re}\langle (\partial_x\cdot\nabla_{\partial_y}-\partial_y\cdot\nabla_{\partial_x})\psi^i,\psi^m\rangle
+\operatorname{Re}\langle\psi^i,(\partial_x\cdot\nabla_{\partial_y}-\partial_y\cdot\nabla_{\partial_x})\psi^m\rangle \\
&-2(\phi^i_x\phi^m_y-\phi^i_y\phi^m_x).
\end{align*}
On the other hand we obtain
\begin{align*}
J^{ij}_x J^{jm}_y-J^{ij}_y J^{jm}_x=&-(\phi^i_x\phi^m_y-\phi^i_y\phi^m_x) \\
&+\operatorname{Re}\langle\psi^i,\partial_x\cdot\psi^j\rangle\operatorname{Re}\langle\psi^j,\partial_y\cdot\psi^m\rangle
-\operatorname{Re}\langle\psi^i,\partial_y\cdot\psi^j\rangle\operatorname{Re}\langle\psi^j,\partial_x\cdot\psi^m\rangle \\
&-\operatorname{Re}\langle\psi^i,\partial_x\cdot\psi^j\rangle\phi^j_y\phi^m
+\operatorname{Re}\langle\psi^j,\partial_y\cdot\psi^m\rangle\phi^i\phi^j_x \\
&+\operatorname{Re}\langle\psi^i,\partial_y\cdot\psi^j\rangle\phi^j_x\phi^m
-\operatorname{Re}\langle\psi^j,\partial_x\cdot\psi^m\rangle\phi^i\phi^j_y.
\end{align*}
Note that all terms proportional to \(\psi^j\phi^j\) drop out since the vector spinors \(\psi\) satisfy \(\psi\perp\phi\).
The result then follows by combining both equations.
\end{proof}
If we also assume that \((\phi,\psi)\colon\Omega\to S^n\) is a smooth \(\kappa\)-Dirac-harmonic map,
we find the following

\begin{Prop}
\label{prop-current-algebra}
Let \((\phi,\psi)\colon\Omega\to S^n\) be a smooth \(\kappa\)-Dirac-harmonic map.
Then the Noether current \(J_\alpha^{im}\) satisfies the following algebra
\begin{align}
\label{noether-algebra-kappa-dh}
\partial_x J_y^{im}&-\partial_y J_x^{im}-2(J^{ij}_x J^{jm}_y-J^{ij}_y J^{jm}_x)=
-\partial_x\langle\psi^i,\partial_y\cdot\psi^m\rangle_{\Sigma M}+\partial_y\langle\psi^i,\partial_x\cdot\psi^m\rangle_{\Sigma M}\\
\nonumber&-2(\operatorname{Re}\langle\psi^i,\partial_x\cdot\psi^j\rangle_{\Sigma M}\operatorname{Re}\langle\psi^j,\partial_y\cdot\psi^m\rangle_{\Sigma M}
-\operatorname{Re}\langle\psi^i,\partial_y\cdot\psi^j\rangle_{\Sigma M}\operatorname{Re}\langle\psi^j,\partial_x\cdot\psi^m\rangle_{\Sigma M})\\
\nonumber&+4\kappa\operatorname{Re}\big(2|\psi|^2\langle \partial_x\cdot \partial_y\cdot\psi^i,\psi^m\rangle_{\Sigma M}
-\langle\psi^j,\psi^i\rangle\langle \partial_x\cdot \partial_y\cdot\psi^j,\psi^m\rangle_{\Sigma M} \\
\nonumber&+\langle\psi^i,\partial_x\cdot \partial_y\cdot\psi^j\rangle\langle\psi^j,\psi^m\rangle_{\Sigma M}\big).
\end{align}
\end{Prop}

\begin{proof}
Multiplying \eqref{kappa-sphere-psi} by \(\partial_x\cdot \partial_y\) we obtain
\begin{align*}
(\partial_x\cdot\nabla_{\partial_y}-\partial_y\cdot\nabla_{\partial_x})\psi^i=&
\partial_y\cdot\psi^j\phi^j_x\otimes\phi^i-\partial_x\cdot\psi^j\phi^j_y\otimes\phi^i \\
&+2\kappa(|\psi|^2\partial_x\cdot \partial_y\cdot\psi^i-\langle\psi^j,\psi^i\rangle \partial_x\cdot \partial_y\cdot\psi^j),
\end{align*}
which yields
\begin{align*}
-&\operatorname{Re}\langle (\partial_x\cdot\nabla_{\partial_y}-\partial_y\cdot\nabla_{\partial_x})\psi^i,\psi^m\rangle
+\operatorname{Re}\langle\psi^i,(\partial_x\cdot\nabla_{\partial_y}-\partial_y\cdot\nabla_{\partial_x})\psi^m\rangle  \\
=&-\operatorname{Re}\langle\partial_y\cdot\psi^j,\psi^m\rangle\phi^j_x\phi^i+\operatorname{Re}\langle\psi^i,\partial_y\cdot\psi^j\rangle\phi^j_x\phi^m 
-\operatorname{Re}\langle\psi^i,\partial_x\cdot\psi^j\rangle\phi^j_y\phi^m+\operatorname{Re}\langle \partial_x\cdot\psi^j,\psi^m\rangle\phi^j_y\phi^i \\
&+2\kappa\operatorname{Re}\big(2|\psi|^2\langle\psi^i,\partial_x\cdot \partial_y\cdot\psi^m\rangle+\langle\psi^j,\psi^i\rangle\langle \partial_x\cdot \partial_y\cdot\psi^j,\psi^m\rangle 
-\langle\psi^i,\partial_x\cdot \partial_y\cdot\psi^j\rangle\langle\psi^m,\psi^j\rangle\big).
\end{align*}
Consequently, the right-hand side of \eqref{noether-current-identity-general} becomes
\begin{align*}
\partial_x& J_y^{im}-\partial_y J_x^{im}-2(J^{ij}_x J^{jm}_y-J^{ij}_y J^{jm}_x)= \\
&-2(\operatorname{Re}\langle\psi^i,\partial_x\cdot\psi^j\rangle\operatorname{Re}\langle\psi^j,\partial_y\cdot\psi^m\rangle
-\operatorname{Re}\langle\psi^i,\partial_y\cdot\psi^j\rangle\operatorname{Re}\langle\psi^j,\partial_x\cdot\psi^m\rangle)\\
&+\operatorname{Re}\langle\partial_y\cdot\psi^j,\psi^m\rangle\phi^j_x\phi^i
-\operatorname{Re}\langle\psi^i,\partial_y\cdot\psi^j\rangle\phi^j_x\phi^m 
+\operatorname{Re}\langle\psi^i,\partial_x\cdot\psi^j\rangle\phi^j_y\phi^m
-\operatorname{Re}\langle \partial_x\cdot\psi^j,\psi^m\rangle\phi^j_y\phi^i\\
&+2\kappa\operatorname{Re}\big(2|\psi|^2\langle \partial_x\cdot \partial_y\cdot\psi^i,\psi^m\rangle-\langle\psi^j,\psi^i\rangle\langle \partial_x\cdot \partial_y\cdot\psi^j,\psi^m\rangle 
+\langle\psi^i,\partial_x\cdot \partial_y\cdot\psi^j\rangle\langle\psi^j,\psi^m\rangle\big).
\end{align*}
Rewriting the terms in the last two lines gives the result.
\end{proof}

\begin{Bem}
Let us make some observations regarding the structure of \eqref{noether-algebra-kappa-dh}.
\begin{enumerate}
\item In the physics literature the Noether algebra \eqref{noether-algebra-kappa-dh} take the simpler form
\begin{align*}
&\partial_x J_y^{im}-\partial_y J_x^{im}+2(J^{ij}_x J^{jm}_y-J^{ij}_y J^{jm}_x)=
-\partial_x\langle\psi^i,\partial_y\cdot\psi^j\rangle_{\Sigma M}+\partial_y\langle\psi^i,\partial_x\cdot\psi^j\rangle_{\Sigma M},
\end{align*}
see for example \cite[p.249]{MR1858211}. To obtain their results physicists make use of so-called \emph{Fierz-identities},
which can be applied to simplify spinorial bilinear terms.
However, physicists usually formulate these identities for Grassmann-valued spinors. 
\item If we think of \(J\) as a connection one-form, then the right hand side of \eqref{noether-algebra-kappa-dh}
gives its curvature.
\end{enumerate}
\end{Bem}

\begin{Bem}
\label{rem-cmc-harmonic-maps}
Suppose that \(\phi\colon\Omega\to S^2\subset\R^3\) is  a harmonic map.
Making use of the conserved Noether current one can show that there exists a map 
\(B\colon\Omega\to\R^3\), unique up to a constant vector, satisfying
\[
\Delta B=2B_x\times B_y.
\]
This equation is well known since it describes a CMC surface when we also require that B is conformal.
More precisely, conformal parametrizations of CMC 1 surfaces are characterized by the system
\begin{equation}
\label{conformal-cmc}
|X_x|^2=|X_y|^2,\qquad\langle X_x,X_y\rangle=0,\qquad \Delta X=2X_x\times X_y. 
\end{equation}
However, we do not know if the map \(B\) is conformal.
If we consider the linear combination \(X_\pm=B\pm u\), for a conformal map \(u\), 
then it turns out that \(X_\pm\) is a solution to the system \eqref{conformal-cmc}.
Assuming that \(u,X_+,X_-\) are immersions, we can associate to a harmonic map
a triple of immersions of surfaces, at a constant distance \(1\) from each other
with \(B(\Omega)\) in the middle (having Gauss curvature \(1\)), and \(X_+(\Omega)\)
and \(X_-(\Omega)\) having mean curvature \(1/2\) at either side.
For more details see \cite[p. 53]{MR1913803} and references therein.
\end{Bem}

According to the last remark the existence of the Noether current for harmonic maps
to spheres leads to a beautiful geometric construction.
In the following we want to discuss if the same holds true for \(\kappa\)-Dirac-harmonic maps
to spheres.

\begin{Prop}
Let \((\phi,\psi)\colon\Omega\to S^n\) be a smooth \(\kappa\)-Dirac-harmonic map.
Then there exist functions \(B^{mi}\) that satisfy
\begin{align}
\label{B-dhmaps}
\Delta B^{mi}=&2(B^{mj}_xB^{ij}_y-B^{mj}_yB^{ji}_x) 
-\partial_x\langle\psi^i,\partial_y\cdot\psi^m\rangle_{\Sigma M}+\partial_y\langle\psi^i,\partial_x\cdot\psi^m\rangle_{\Sigma M}\\
\nonumber&-2(\operatorname{Re}\langle\psi^i,\partial_x\cdot\psi^j\rangle_{\Sigma M}\operatorname{Re}\langle\psi^j,\partial_y\cdot\psi^m\rangle_{\Sigma M}
-\operatorname{Re}\langle\psi^i,\partial_y\cdot\psi^j\rangle_{\Sigma M}\operatorname{Re}\langle\psi^j,\partial_x\cdot\psi^m\rangle_{\Sigma M})\\
\nonumber&+4\kappa\operatorname{Re}\big(2|\psi|^2\langle \partial_x\cdot \partial_y\cdot\psi^i,\psi^m\rangle_{\Sigma M}
-\langle\psi^j,\psi^i\rangle\langle \partial_x\cdot \partial_y\cdot\psi^j,\psi^m\rangle_{\Sigma M} \\
\nonumber&+\langle\psi^i,\partial_x\cdot \partial_y\cdot\psi^j\rangle\langle\psi^j,\psi^m\rangle_{\Sigma M}\big)
\end{align}
and
\begin{align*}
|B_x|^2=&|\langle\psi^i,\partial_y\cdot\psi^m\rangle_{\Sigma M}|^2-2|\phi_y|^2,\qquad |B_y|^2=|\langle\psi^i,\partial_x\cdot\psi^m\rangle_{\Sigma M}|^2-2|\phi_x|^2, \\
&\langle B_x,B_y\rangle=\langle\psi^i,\partial_x\cdot\psi^m\rangle_{\Sigma M}\langle\psi^i,\partial_y\cdot\psi^m\rangle_{\Sigma M}-2\phi^i_y\phi^i_x.
\end{align*}

\end{Prop}
\begin{proof}
Since the Noether current \eqref{noether-current-kappa-sphere} is divergence-free, we have
\begin{align*}
\partial_x\big(\langle\psi^i,\partial_x\cdot\psi^m\rangle-(\phi^i_x\phi^m-\phi^i\phi^m_x)\big)
+\partial_y\big(\langle\psi^i,\partial_y\cdot\psi^m\rangle-(\phi^i_y\phi^m-\phi^i\phi^m_y)\big)=0.
\end{align*}
Hence, there must exist functions \(B^{mi}\) that satisfy 
\begin{align*}
B^{mi}_x=&\langle\psi^i,\partial_y\cdot\psi^m\rangle-\phi^i_y\phi^m+\phi^i\phi^m_y, \\
B^{mi}_y=&-\langle\psi^i,\partial_x\cdot\psi^m\rangle+\phi^i_x\phi^m-\phi^i\phi^m_x.
\end{align*}
By a direct calculation we find
\begin{align*}
\Delta B^{mi}=\partial_x J^{im}_y-\partial_y J^{im}_x
\end{align*}
and the result follows by Proposition \ref{prop-current-algebra}.

\end{proof}

\begin{Bem}
It is obvious that we do not get a nice geometric configuration as for harmonic maps to spheres from \eqref{B-dhmaps}
due to the presence of the spinors.
\end{Bem}

Although the regularity theory for Dirac-harmonic maps with curvature term is fully developed by now \cite{MR3333092}
we want to point out how the Noether current can be used to establish the continuity of the map \(\phi\),
whenever we are given a weak solutions of \eqref{kappa-critical-phi}, \eqref{kappa-critical-psi} with a spherical target.

For harmonic maps to spheres this was first noted in \cite[Proposition 2.1]{MR1078114},
and for Dirac-harmonic maps to spheres in \cite[Remark 4.4]{MR2176464} without referring to the Noether current.
\begin{Prop}
Let \((\phi,\psi)\colon\Omega\to S^n\) be a weak \(\kappa\)-Dirac-harmonic map.
There exists \(M\in W^{1,2}(\Omega,\R^{(n+1)\times(n+1)})\) such that
\begin{align}
-\Delta\phi=\frac{\partial\phi}{\partial x}\frac{\partial M}{\partial y}-\frac{\partial\phi}{\partial y}\frac{\partial M}{\partial x}
\end{align}
holds.
\end{Prop}
\begin{proof}
We calculate (in a distributional sense)
\begin{align*}
\Delta\phi^m=&-|d\phi|^2\phi^m-\operatorname{Re}\langle\psi^m,e_\alpha\cdot\psi^i\rangle\phi^i_\alpha \\
=&-(\phi^i_\alpha\phi^m-\phi^i\phi^m_\alpha)\phi^i_\alpha-\operatorname{Re}\langle\psi^m,e_\alpha\cdot\psi^i\rangle\phi^i_\alpha \\
=&-J^{im}_\alpha\phi^i_\alpha.
\end{align*}
Since the Noether current \(J_\alpha^{im}\) is conserved there exist functions 
\(M^{mi}\) on \(\Omega\) satisfying
\[
J_x^{im}=\frac{\partial M^{im}}{\partial y} ,\qquad J_y^{im}=-\frac{\partial M^{im}}{\partial x},
\]
which completes the proof.
\end{proof}

\begin{Cor}
This yields continuity of \(\phi\) via Wente's Lemma \cite{MR0243467} for all values of \(\kappa\).
\end{Cor}

\begin{Bem}
If one also considers a two-form contribution in the action functional as in \cite{MR3305429}
then the Noether current is no longer conserved.
This is not surprising from a physical point of view: The two-form potential
in the energy functional is used to model a (generalized) external magnetic field.
However, a magnetic field always destroys the rotational symmetry of a system since
it introduces a preferred direction.
\end{Bem}

\begin{Bem}
The norm of the Noether current satisfies
\begin{equation}
|J^{im}_\alpha|^2=|\operatorname{Re}\langle\psi^i,e_\alpha\cdot\psi^m\rangle_{\Sigma M}|^2+|\phi_\alpha|^2.
\end{equation}
\end{Bem}
\begin{proof}
We calculate
\begin{align*}
|J^{im}_\alpha|^2=|\operatorname{Re}\langle\psi^i,e_\alpha\cdot\psi^m\rangle|^2+2|\phi_\alpha|^2
-2(\operatorname{Re}\langle\psi^i,e_\alpha\cdot\psi^m\rangle)(\phi^i_\alpha\phi^m-\phi^i\phi^m_\alpha).
\end{align*}
Note that the mixed terms vanish since \(\psi\perp\phi\).
\end{proof}

In the following we want to explore the limit \(\phi=const\), which 
is well-known in the physics literature as the Gross-Neveu model.

\subsection{The Gross Neveu model and CMC surfaces}
The Gross-Neveu model \cite{PhysRevD.10.3235} is a model for interacting massive fermions in two dimensions.
For its mathematical study let \((M,h)\) be a closed Riemannian spin surface.
For a geometric treatment of the Gross-Neveu model on complete Riemannian
manifolds see \cite{branding2016nonlinear}.

Its energy functional is given by
\begin{align}
\label{energy-gross-neveu}
E(\psi)=\int_M(\langle\psi,\p\psi\rangle-\lambda|\psi|^2-\frac{\kappa}{2}|\psi|^4\rangle)dM,
\end{align}
where \(\lambda\) and \(\kappa\) are real-valued parameters and
\(\psi\in\Gamma(\Sigma M\otimes\R^{q})\).

The critical points of \eqref{energy-gross-neveu} are given by
\begin{align}
\label{gross-neveu-psi}
\p\psi^i=\lambda\psi^i+\kappa|\psi|^2\psi^i.
\end{align}

The analytic aspects of such kind of nonlinear Dirac equations
have been studied in \cite{MR2661574, MR2390834}.

\begin{Lem}
Let \(\psi\in\Gamma(\Sigma M\otimes\R^{q})\) be a solution of \eqref{gross-neveu-psi}.
Then the Noether current 
\begin{align*}
\label{noether-currrent-gross-neveu}
J^{im}_\alpha:=\langle\psi^i,e_\alpha\cdot\psi^m\rangle_{\Sigma M}
\end{align*}
is conserved, that is 
\begin{align}
\nabla_{e_\alpha}J^{im}_\alpha=0.
\end{align}
\end{Lem}
\begin{proof}
This follows by a direct calculation.
\end{proof}

In order to derive the corresponding Noether algebra we need an algebraic relation
for the spinorial bilinear terms. Since we are only interested in a local statement,
we choose a local trivialization of the spinor bundle \(\Sigma M\) such that we 
can work with complex-valued functions.

\begin{Lem}
Let \(\psi^i,\psi^j,\psi^m\in\Gamma(\Sigma M)\).
Then the following algebraic identity holds
\begin{align}
\label{fierz}\nonumber
\langle\psi^i,\partial_x\cdot\psi^j\rangle_{\Sigma M}&\langle\psi^j,\partial_y\cdot\psi^m\rangle_{\Sigma M}
-\langle\psi^i,\partial_y\cdot\psi^j\rangle_{\Sigma M}\langle\psi^j,\partial_x\cdot\psi^m\rangle_{\Sigma M}
=2\langle\psi^i,\partial_x\cdot\partial_y\cdot\psi^m\rangle_{\Sigma M}|\psi|_{\Sigma M}^2 \\
&+|P_-\psi^j|_{\Sigma M}^2\langle P_-\psi^i,P_-\psi^m\rangle_{\Sigma M}
-|P_+\psi^j|^2\langle P_+\psi^i,P_+\psi^m\rangle_{\Sigma M},
\end{align}
where \(P_\pm\) is defined in \eqref{projector-chiralities}.
\end{Lem}

\begin{proof}
We prove the identity via a local calculation. 
Locally, the spinors \(\psi^i, i=1,\ldots,q\) can be thought of as \(\C^2\)-valued functions and we choose
\begin{align*}
\psi^i=
\begin{pmatrix}
a_1  \\
a_2
\end{pmatrix},
\qquad
\psi^j=
\begin{pmatrix}
b_1  \\
b_2
\end{pmatrix},
\qquad
\psi^m=
\begin{pmatrix}
c_1  \\
c_2
\end{pmatrix},
\end{align*}
where \(a_1,a_2,b_1,b_2,c_1,c_2\) are complex-valued functions.
In addition, Clifford multiplication with \(\partial_x\) and \(\partial_y\) can be expressed as
multiplication with the matrices
\begin{align*}
\partial_x\cdot=
\begin{pmatrix}
0 & 1 \\
-1 & 0
\end{pmatrix},
\qquad
\partial_y\cdot=
\begin{pmatrix}
0 & i \\
i & 0
\end{pmatrix}.
\end{align*}

By a direct calculation using the standard hermitian scalar product on \(\C^2\) we find
\begin{align*}
\langle\psi^i,\partial_x\cdot\psi^j\rangle_{\Sigma M}\langle\psi^j,\partial_y\cdot\psi^m\rangle_{\Sigma M}
-\langle\psi^i,\partial_y\cdot\psi^j\rangle_{\Sigma M}\langle\psi^j,\partial_x\cdot\psi^m\rangle_{\Sigma M}
=2i(a_2\bar{c}_2b_1\bar{b}_1-a_1\bar{c}_1b_2\bar{b}_2)
\end{align*}
and also
\begin{align*}
2\langle\psi^i,\partial_x\cdot\partial_y\cdot\psi^m\rangle_{\Sigma M}|\psi|_{\Sigma M}^2
=2i(-a_1\bar{c}_1b_1\bar{b}_1+a_2\bar{c}_2b_1\bar{b}_1-a_1\bar{c}_1b_2\bar{b}_2+a_2\bar{c}_2b_2\bar{b}_2).
\end{align*}
We require that
\[
a_1\bar{c}_1b_1\bar{b}_1=a_2\bar{c}_2b_2\bar{b}_2,
\]
which is equivalent to
\[
|P_-\psi^j|_{\Sigma M}^2\langle P_-\psi^i,P_-\psi^m\rangle_{\Sigma M}=|P_+\psi^j|_{\Sigma M}^2\langle P_+\psi^i,P_+\psi^m\rangle_{\Sigma M},
\]
completing the proof.
\end{proof}

\begin{Cor}
Let \(\psi^i,\psi^j,\psi^m\in\Gamma(\Sigma M)\). If in addition
\begin{align}
\label{majorana-condition}
|P_-\psi^j|_{\Sigma M}^2\langle P_-\psi^i,P_-\psi^m\rangle_{\Sigma M}=|P_+\psi^j|_{\Sigma M}^2\langle P_+\psi^i,P_+\psi^m\rangle_{\Sigma M}
\end{align}
holds, then \eqref{fierz} simplifies to
\begin{align}
\label{fierz-majorana}
\langle\psi^i,\partial_x\cdot\psi^j\rangle_{\Sigma M}\langle\psi^j,\partial_y\cdot\psi^m\rangle_{\Sigma M}
-\langle\psi^i,\partial_y\cdot\psi^j\rangle_{\Sigma M}\langle\psi^j,\partial_x\cdot\psi^m\rangle_{\Sigma M} 
=2\langle\psi^i,\partial_x\cdot\partial_y\cdot\psi^m\rangle_{\Sigma M}|\psi|^2.
\end{align}
\end{Cor}

\begin{Lem}
Let \(\psi\) be a smooth solution of \eqref{gross-neveu-psi}.
Moreover, suppose that \eqref{majorana-condition} holds.
Then the Noether current \eqref{noether-currrent-gross-neveu} satisfies the following algebra
\begin{align}
\label{noether-algebra-gross-neveu}
\partial_xJ_y^{im}-\partial_yJ_x^{im}-\kappa(J^{ij}_xJ^{jm}_y-J^{ij}_yJ^{jm}_x)=2\lambda\langle\psi^i,\partial_x\cdot\partial_y\cdot\psi^m\rangle_{\Sigma M}.
\end{align}
\end{Lem}
\begin{proof}
By a direct calculation we find
\begin{align*}
\partial_xJ^{im}_y-\partial_yJ^{im}_x
=&\langle(\partial_x\cdot\nabla_{\partial_y}-\partial_y\cdot\nabla_{\partial_x})\psi^i,\psi^m\rangle
+\langle\psi^i,(\partial_x\cdot\nabla_{\partial_y}-\partial_y\cdot\nabla_{\partial_x})\psi^m\rangle \\
=&2\kappa|\psi|^2\langle\psi^i,\partial_x\cdot\partial_y\cdot\psi^m\rangle
+2\lambda\langle\psi^i,\partial_x\cdot\partial_y\cdot\psi^m\rangle.
\end{align*}
In addition, we find
\begin{align*}
J^{ij}_xJ^{jm}_y-J^{ij}_yJ^{jm}_x=\langle\psi^i,\partial_x\cdot\psi^j\rangle\langle\psi^j,\partial_y\cdot\psi^m\rangle
-\langle\psi^i,\partial_y\cdot\psi^j\rangle\langle\psi^j,\partial_x\cdot\psi^m\rangle.
\end{align*}
The claim then follows from the Fierz identity \eqref{fierz-majorana}.
\end{proof}

\begin{Bem}
It is obvious that the Noether algebra has the form of a zero-curvature condition for \(J\) when we are considering
the massless Gross-Neveu-model, that is \(\lambda=0\).
\end{Bem}

\begin{Prop}
Let \(\psi\) be a smooth solution of \eqref{gross-neveu-psi}.
In addition, suppose that \eqref{majorana-condition} holds.
Then there exist functions \(B\) that satisfy
\begin{align}
\label{b-gross-neveu} 
\Delta B^{mi}=\kappa(B^{mj}_xB^{ji}_y-B^{mj}_yB^{ji}_x)+2\lambda\langle\psi^i,\partial_x\cdot\partial_y\cdot\psi^m\rangle_{\Sigma M}.
\end{align}
\end{Prop}
\begin{proof}
Since the Noether current is divergence free, there exist functions \(B\) that satisfy
\begin{align*}
-B^{ij}_y=\langle\psi^i,\partial_x\cdot\psi^j\rangle,\qquad B^{ij}_x=\langle\psi^j,\partial_y\cdot\psi^j\rangle.
\end{align*}
Thus, a direct calculation yields
\[
\Delta B=\partial_xJ_y-\partial_yJ_x
\]
and the result follows from \eqref{noether-algebra-gross-neveu}.
\end{proof}

\begin{Bem}
In the case that \(\lambda=0\) the Noether algebra \eqref{b-gross-neveu} satisfies
the equation for a CMC surface.
However, \(B\) is not conformal, since
\[
|B_x|^2\neq |B_y|^2,\qquad \langle B_x,B_y\rangle\neq 0
\]
without posing further assumptions.
\end{Bem}

We can again use the Noether current to establish some regularity result.
However, the regularity of weak solutions of \eqref{gross-neveu-psi} is already well-understood,
see \cite{MR2550205,MR2661574}.

\begin{Bem}
Let \(\psi\in W^{1,\frac{4}{3}}(\Sigma M)\) be a distributional solution of \eqref{gross-neveu-psi}
with \(\lambda=0\).
Again, by application of the Wente Lemma we find that the map \(B\) is continuous
since
\[
\|B_x\|_{L^2}\leq C\|\psi\|^2_{L^4}\leq C\|\nabla\psi\|^2_{W^{1,\frac{4}{3}}},
\]
where the last estimate follows from the Sobolev embedding in two-dimensions.
However, we cannot use the statement on the regularity of \(B\) to gain regularity for \(\psi\).
\end{Bem}

\begin{Bem}
It is obvious that the algebra of the Noether current for harmonic maps to spheres and
for the massless \((\lambda=0)\) Gross-Neveu model \eqref{noether-algebra-gross-neveu} is the same. 
This fact suggests that both models describe similar geometric and physical phenomena.

In physics this fact is often referred to as \emph{bosonization}, which reflects the fact that
a combination of two fermions behaves like a boson.

In geometric terms we have seen the relationship between harmonic maps to spheres and CMC
surfaces in Remark \ref{rem-cmc-harmonic-maps}.
On the other hand, it is also well-known that the solutions of 
nonlinear Dirac-equations of the form \eqref{gross-neveu-psi} with \(\lambda=0\) describe CMC surfaces
from the universal covering of \(M\) in \(\R^3\), see \cite{MR1653146,MR2550205}.
More precisely, we have the following bijection
\begin{align*}
\{\text{Solutions of } \p\psi=H|\psi|^2\psi\}/\pm 1 \leftrightarrow
\{
\text{Conformal periodic H-immersions }\tilde M\subset\R^3 \\
\text{with branching points of even order}
\},
\end{align*}
where \(\tilde M\) denotes the universal covering of \(M\).
\end{Bem}

\section{Conservation laws for targets with Killing vector fields}
In this section we discuss how to generalize the notion of the Noether current to
target manifolds that possess Killing vector fields. 
The approach that we take here is different compared to the one that is usually 
taken in the physics literature. We derive the Noether current by assuming that the target
manifold admits Killing vector fields, whereas in the physics literature the Noether current is
obtained by considering symmetries acting on the fields that leave the action invariant.

Let \(\xi\) be a diffeomorphism that generates a one-parameter family of vector fields \(X\).
Then we know that
\begin{align*}
\frac{d}{dt}\big|_{t=0}\xi^\ast g=\cL_Xg.
\end{align*}

This enables us to give the following 
\begin{Dfn}
Let \(\xi\) be a diffeomorphism that generates a one-parameter family of vector fields \(X\).
Then we say that \(X\) generates a symmetry for the action \(E_\kappa(\phi,\psi,\xi^\ast g)\) if
\begin{align*}
\frac{d}{dt}\big|_{t=0}E_\kappa(\phi,\psi,\xi^\ast g)=\int_M\cL_X(|d\phi|^2+\langle\psi,\D\psi\rangle+\kappa\langle R^N(\psi,\psi)\psi,\psi\rangle)dM=0,
\end{align*}
where the Lie-derivative is acting on the metric \(g\).
\end{Dfn}

Note that if \(X\) generates an isometry then \(\cL_Xg=0\) such that we have to require the existence of Killing
vector fields on the target.

\begin{Lem}
Let \((\phi,\psi)\) be a smooth \(\kappa\)-Dirac-harmonic map to a target with Killing vector fields.
Then the Lie-derivative acting on the metric \(g\) of the energy density is given by
\begin{align}
\label{lie-derivative-energy}
\cL_X(|d\phi|^2+\langle\psi,\D\psi\rangle+\kappa\langle R^N(\psi,\psi)\psi,\psi\rangle)
\nonumber=&2\nabla_{e_\alpha}\langle d\phi(e_\alpha),X(\phi)\rangle-\langle R^N(e_\alpha\cdot\psi,\psi)d\phi(e_\alpha),X(\phi)\rangle \\
&+2\kappa\langle\psi^i,\psi^k\rangle\langle\psi^j,\psi^l\rangle(R_{ijsl}\nabla_kX_s-R_{sjkl}\nabla_sX_i).
\end{align}
\end{Lem}
\begin{proof}
We calculate (with \(x_\alpha\) being local coordinates on \(M\))
\begin{align*}
\cL_X|d\phi|^2=&(\cL_Xg)_{ij}\frac{\partial\phi^i}{\partial x^\alpha}\frac{\partial\phi^j}{\partial x^\beta}h^{\alpha\beta} \\
=&2\nabla_iX_j\frac{\partial\phi^i}{\partial x^\alpha}\frac{\partial\phi^j}{\partial x^\beta}h^{\alpha\beta} \\
=&2\langle d\phi(e_\alpha),\nabla_{e_\alpha}(X(\phi))\rangle \\
=&2\nabla_{e_\alpha}\langle d\phi(e_\alpha),X(\phi)\rangle-2\langle\tau(\phi),X(\phi)\rangle \\
=&2\nabla_{e_\alpha}\langle d\phi(e_\alpha),X(\phi)\rangle-\langle R^N(e_\alpha\cdot\psi,\psi)d\phi(e_\alpha),X(\phi)\rangle
-\kappa\langle(\nabla_X R)(\psi,\psi)\psi,\psi\rangle,
\end{align*}
where we used that \(\phi\) is a solution of \eqref{kappa-critical-phi} in the last step.

As a second step we calculate
\begin{align*}
\cL_X\langle\psi,\D\psi\rangle=&\langle\psi^i,(\D\psi)^j\rangle(\cL_Xg)_{ij} \\
=&\langle\psi^i,(\D\psi)^j\rangle(\nabla_iX_j+\nabla_jX_i)\\
=&-2\kappa R_{abcj}\langle\psi^a,\psi^c\rangle\langle\psi^i,\psi^b\rangle\rangle(\nabla_iX_j+\nabla_jX_i),
\end{align*}
where we used that \(\psi\) is a solution of \eqref{kappa-critical-psi}.
Recall the formula for the Lie-derivative of the Riemann curvature tensor
\begin{align*}
\cL_XR_{ijkl}=(\nabla_XR)_{ijkl}+R_{sjkl}\nabla_iX_s+R_{iskl}\nabla_jX_s+R_{ijsl}\nabla_kX_s+R_{ijks}\nabla_lX_s.
\end{align*}
Consequently, we find
\begin{align*}
(\cL_XR_{ijkl})\langle\psi^i,\psi^k\rangle\langle\psi^j,\psi^l\rangle
=&(\nabla_XR_{ijkl})\langle\psi^i,\psi^k\rangle\langle\psi^j,\psi^l\rangle \\
&+2\langle\psi^i,\psi^k\rangle\langle\psi^j,\psi^l\rangle(R_{sjkl}\nabla_iX_s+R_{ijsl}\nabla_kX_s).
\end{align*}
Adding up the three contributions yields the result.
\end{proof}

Note that the first term on the right hand side of \eqref{lie-derivative-energy} already is in divergence form,
which is what we need to obtain a conservation law. To rewrite the other terms on the right hand side of \eqref{lie-derivative-energy} 
we need the following
\begin{Lem}
Let \(X\) be a Killing vector field on a Riemannian manifold \((N,g)\),
then the following formula holds
\begin{align}
\label{killing-curvature-tensor-identity}
-R^N(X,Y)Z=\nabla^2_{Y,Z}X.
\end{align}
\end{Lem}
\begin{proof}
A proof can be found in \cite[p.242, Lemma 33]{MR2243772}.
\end{proof}

From now on \(X\) will always denote a Killing vector field on \(N\).

First, we will give a conservation law for Dirac-harmonic maps, that is for solutions of \eqref{kappa-critical-phi}, \eqref{kappa-critical-psi}
with \(\kappa=0\).

\begin{Satz}
Let \((\phi,\psi)\) be a smooth Dirac-harmonic map to a target with Killing vector fields.
Then the current defined by
\begin{align}
\label{noether-dhmap-general}
J_\alpha:=2\langle d\phi(e_\alpha),X(\phi)\rangle_{\phi^\ast TN}-\langle\nabla_\psi X(\phi),e_\alpha\cdot\psi\rangle_{\Sigma M\otimes\phi^\ast TN}
\end{align}
is conserved, that is \(\nabla_{e_\alpha} J_\alpha=0\). Here, the notation is to be understood as
\[
\langle\nabla_\psi X(\phi),e_\alpha\cdot\psi\rangle:=\langle\psi^i,e_\alpha\cdot\psi^j\rangle\nabla_{\partial_{y^i}}X_j,
\]
where \(\partial_{y^i}\) is a local basis of \(TN\).
\end{Satz}
\begin{proof}
By a direct calculation we find
\begin{align*}
\nabla_{e_\alpha}\langle d\phi(e_\alpha),X(\phi)\rangle=&2\langle\tau(\phi),X(\phi)\rangle+\underbrace{\langle d\phi(e_\alpha),\nabla_{d\phi(e_\alpha)}X(\phi)\rangle}_{=0} 
=\langle R^N(e_\alpha\cdot\psi,\psi)d\phi(e_\alpha),X(\phi)\rangle.
\end{align*}
On the other hand we get
\begin{align*}
\nabla_{e_\alpha}\langle\nabla_\psi X(\phi),e_\alpha\cdot\psi\rangle=&\langle\nabla^2_{d\phi(e_\alpha),\psi} X(\phi),e_\alpha\cdot\psi\rangle
+\langle\nabla_\psi X(\phi),\underbrace{\D\psi}_{=0}\rangle \\
=&\langle R^N(e_\alpha\cdot\psi,\psi)d\phi(e_\alpha),X(\phi)\rangle,
\end{align*}
where we applied \eqref{killing-curvature-tensor-identity} in the last step. The result then follows by combining the two equations.
\end{proof}

As for harmonic maps \cite{MR1085633}, we can use the conserved current \(J_\alpha\) to study the regularity of weak Dirac-harmonic maps.

\begin{Prop}
Suppose there exists a finite dimensional Lie group \(G\) which acts transitively
on \(N\) by isometries. Then for a given weak Dirac-harmonic map we can deduce that \(\phi\) is continuous.
\end{Prop}
\begin{proof}
This follows directly as in \cite[Theorem A]{MR1085633}.
\end{proof}

\begin{Bem}
Let \((\phi,\psi)\) be a smooth Dirac-harmonic map from a surface to a target with Killing vector fields.
Then a direct calculation yields
\begin{align*}
\nabla_{\partial_y}J_x-\nabla_{\partial_x}J_y
=&2\langle\frac{\partial\phi}{\partial x},\nabla_{d\phi(\partial_y)}X\rangle-2\langle\frac{\partial\phi}{\partial y},\nabla_{d\phi(\partial_x)}X\rangle \\
&+\langle R^N(\partial_x\cdot\psi,\psi)\frac{\partial\phi}{\partial y}-R^N(\partial_y\cdot\psi,\psi)\frac{\partial\phi}{\partial x},X\rangle.
\end{align*}
Due to the non-trivial curvature of the target manifold we cannot rewrite the right hand side of this equation
in terms of the current \(J_\alpha\) as we could do in the case of a spherical target.
\end{Bem}

Finally, we give a conservation law for solutions of \eqref{kappa-critical-phi}, \eqref{kappa-critical-psi}
in the case of \(\kappa\neq 0\). It turns out, that we have to impose additional restrictions on the curvature
of the target manifold.

\begin{Satz}
Let \((\phi,\psi)\) be a smooth \(\kappa\)-Dirac-harmonic map to a target with Killing vector fields.
Then the current defined by
\begin{align}
J_\alpha:=2\langle d\phi(e_\alpha),X(\phi)\rangle-\langle\nabla_\psi X(\phi),e_\alpha\cdot\psi\rangle
\end{align} 
is conserved, if \(N\) has constant curvature \(K\).
\end{Satz}
\begin{proof}
Performing a similar calculation as before, we find
\begin{align*}
\nabla_{e_\alpha}J_\alpha=-\kappa\langle(\nabla_X R^N)(\psi,\psi)\psi,\psi\rangle
+2\kappa\langle\nabla_\psi X(\phi),R^N(\psi,\psi)\psi\rangle.
\end{align*}
In general, we cannot expect to rewrite the right hand side as a divergence term 
since the right hand side of \eqref{lie-derivative-energy} does not vanish for \(\kappa\neq 0\).
However, the first term on the right hand side vanishes by assumption. For the second term we rewrite
\begin{align*}
\kappa\langle\nabla_\psi X,R^N(\psi,\psi)\psi\rangle&=
R_{ijkl}\nabla_sX_i\langle\psi^j,\psi^l\rangle\langle\psi^k,\psi^s\rangle \\
&=K\nabla_sX_i(\langle\psi^j,\psi^i\rangle\langle\psi^j,\psi^s\rangle-|\psi|^2\langle\psi^i,\psi^s\rangle)\\
&=0,
\end{align*}
where we used the assumption that \(N\) has constant curvature \(K\) in the second step.
The above expression vanishes due to the skew-symmetry of \(\nabla_sX_i\).
\end{proof}

\begin{Bem}
On a closed Riemannian surface we can always find a metric of constant curvature \(K\)
due to the uniformization theorem. Consequently, the last Theorem always holds for 
smooth \(\kappa\)-Dirac-harmonic maps to a closed two-dimensional target.
\end{Bem}

\par\medskip
\textbf{Acknowledgements:}
The author gratefully acknowledges the support of the Austrian Science Fund (FWF) 
through the START-Project Y963-N35 of Michael Eichmair.
\bibliographystyle{plain}
\bibliography{mybib}

\begin{thebibliography}{10}

\bibitem{MR1858211}
Elcio Abdalla, M.~Cristina~B. Abdalla, and Klaus~D. Rothe.
\newblock {\em Non-perturbative methods in 2 dimensional quantum field theory}.
\newblock World Scientific Publishing Co., Inc., River Edge, NJ, second
  edition, 2001.

\bibitem{MR626710}
Luis Alvarez-Gaum\'e and Daniel~Z. Freedman.
\newblock Geometrical structure and ultraviolet finiteness in the
  supersymmetric {$\sigma $}-model.
\newblock {\em Comm. Math. Phys.}, 80(3):443--451, 1981.

\bibitem{MR2550205}
Bernd Ammann.
\newblock The smallest {D}irac eigenvalue in a spin-conformal class and cmc
  immersions.
\newblock {\em Comm. Anal. Geom.}, 17(3):429--479, 2009.

\bibitem{MR3070562}
Bernd Ammann and Nicolas Ginoux.
\newblock Dirac-harmonic maps from index theory.
\newblock {\em Calc. Var. Partial Differential Equations}, 47(3-4):739--762,
  2013.

\bibitem{MR1224089}
Christian B\"ar.
\newblock Real {K}illing spinors and holonomy.
\newblock {\em Comm. Math. Phys.}, 154(3):509--521, 1993.

\bibitem{MR3305429}
Volker Branding.
\newblock Magnetic {D}irac-harmonic maps.
\newblock {\em Anal. Math. Phys.}, 5(1):23--37, 2015.

\bibitem{MR3333092}
Volker Branding.
\newblock Some aspects of {D}irac-harmonic maps with curvature term.
\newblock {\em Differential Geom. Appl.}, 40:1--13, 2015.

\bibitem{MR3493217}
Volker Branding.
\newblock Dirac-harmonic maps with torsion.
\newblock {\em Commun. Contemp. Math.}, 18(4):1550064, 19, 2016.

\bibitem{branding2016nonlinear}
Volker Branding.
\newblock Nonlinear {D}irac equations, monotonicity formulas and {L}iouville
  theorems.
\newblock {\em arXiv preprint arXiv:1605.03453}, 2016.

\bibitem{MR2370260}
Q.~Chen, J.~Jost, and G.~Wang.
\newblock Liouville theorems for {D}irac-harmonic maps.
\newblock {\em J. Math. Phys.}, 48(11):113517, 13, 2007.

\bibitem{MR2176464}
Qun Chen, J\"urgen Jost, Jiayu Li, and Guofang Wang.
\newblock Regularity theorems and energy identities for {D}irac-harmonic maps.
\newblock {\em Math. Z.}, 251(1):61--84, 2005.

\bibitem{MR2262709}
Qun Chen, J\"urgen Jost, Jiayu Li, and Guofang Wang.
\newblock Dirac-harmonic maps.
\newblock {\em Math. Z.}, 254(2):409--432, 2006.

\bibitem{MR2390834}
Qun Chen, J\"urgen Jost, and Guofang Wang.
\newblock Nonlinear {D}irac equations on {R}iemann surfaces.
\newblock {\em Ann. Global Anal. Geom.}, 33(3):253--270, 2008.

\bibitem{MR1653146}
Thomas Friedrich.
\newblock On the spinor representation of surfaces in {E}uclidean {$3$}-space.
\newblock {\em J. Geom. Phys.}, 28(1-2):143--157, 1998.

\bibitem{PhysRevD.10.3235}
David~J. Gross and Andr\'e Neveu.
\newblock Dynamical symmetry breaking in asymptotically free field theories.
\newblock {\em Phys. Rev. D}, 10:3235--3253, Nov 1974.

\bibitem{MR1078114}
Fr\'ed\'eric H\'elein.
\newblock R\'egularit\'e des applications faiblement harmoniques entre une
  surface et une sph\`ere.
\newblock {\em C. R. Acad. Sci. Paris S\'er. I Math.}, 311(9):519--524, 1990.

\bibitem{MR1085633}
Fr\'ed\'eric H\'elein.
\newblock Regularity of weakly harmonic maps from a surface into a manifold
  with symmetries.
\newblock {\em Manuscripta Math.}, 70(2):203--218, 1991.

\bibitem{MR2389639}
Fr\'ed\'eric H\'elein and John~C. Wood.
\newblock Harmonic maps.
\newblock In {\em Handbook of global analysis}, pages 417--491, 1213. Elsevier
  Sci. B. V., Amsterdam, 2008.

\bibitem{MR1913803}
Fr\'ed\'eric~and H\'elein.
\newblock {\em Harmonic maps, conservation laws and moving frames}, volume 150
  of {\em Cambridge Tracts in Mathematics}.
\newblock Cambridge University Press, Cambridge, second edition, 2002.
\newblock Translated from the 1996 French original, With a foreword by James
  Eells.

\bibitem{MR2316716}
Idrisse Khemar.
\newblock Supersymmetric harmonic maps into symmetric spaces.
\newblock {\em J. Geom. Phys.}, 57(8):1601--1630, 2007.

\bibitem{MR1031992}
H.~Blaine Lawson, Jr. and Marie-Louise Michelsohn.
\newblock {\em Spin geometry}, volume~38 of {\em Princeton Mathematical
  Series}.
\newblock Princeton University Press, Princeton, NJ, 1989.

\bibitem{MR2243772}
Peter Petersen.
\newblock {\em Riemannian geometry}, volume 171 of {\em Graduate Texts in
  Mathematics}.
\newblock Springer, New York, second edition, 2006.

\bibitem{MR2661574}
Changyou Wang.
\newblock A remark on nonlinear {D}irac equations.
\newblock {\em Proc. Amer. Math. Soc.}, 138(10):3753--3758, 2010.

\bibitem{MR2544729}
Changyou Wang and Deliang Xu.
\newblock Regularity of {D}irac-harmonic maps.
\newblock {\em Int. Math. Res. Not. IMRN}, (20):3759--3792, 2009.

\bibitem{MR0243467}
Henry~C. Wente.
\newblock An existence theorem for surfaces of constant mean curvature.
\newblock {\em J. Math. Anal. Appl.}, 26:318--344, 1969.

\bibitem{MR1264181}
J.~C. Wood.
\newblock Harmonic maps into symmetric spaces and integrable systems.
\newblock In {\em Harmonic maps and integrable systems}, Aspects Math., E23,
  pages 29--55. Friedr. Vieweg, Braunschweig, 1994.

\end{thebibliography}
\end{document}